# Numerical Solution of Multiple Nonlinear Volterra Integral Equations

S. A. Belbas(*), Mathematics Department, University of Alabama, Tuscaloosa, AL 35487-0350, USA. *e-mail*: SBELBAS@AS.UA.EDU

Yuriy Bulka, Department of Mathematics, Austin Peay State University, P. O. Box 4626, Clarksville, TN 37044, USA. *e-mail:* BULKAY@APSU.EDU

(*) Corresponding author.

Abstract. We analyze a discretization method for solving nonlinear integral equations that contain multiple integrals. These equations include integral equations with a Volterra series, instead of a single integral term, on one side of the equation. We prove existence and uniqueness of solutions, and convergence and estimates of the order of convergence for the numerical methods of solution.

# 1. INTRODUCTION

In this work we briefly outline some analytical results and then investigate in detail a numerical method for solving multiple nonlinear Volterra integral equations. These are extensions of one-dimensional nonlinear Volterra integral equation, i.e. the equation of the form

$$x(t) = x_0(t) + \int_0^t f(t,s,x(s))ds \qquad (1.1)$$

Eq. (1.1) is widely used to model causal systems with memory arising in different fields. Both analytical properties and numerical solution of single Volterra integral equation (1.1) has been studied, see for example Corduneanu [C1]. However, single Volterra equations do not represent the most general form of causal systems with memory. In recent years, there has been growing interest in functional equations with general causal memory terms [C2]. In contrast to abstract functional equations, in this paper we deal with a concrete instance of causal equations with memory that are more general than the usual Volterra equations, namely equations with multiple integral terms. It is known (see, for example, Rugh [R], Schetzen [S]) that in the context of general causal systems with memory polynomial Volterra operators of the type

$$V(x)(t) = x_0(t) + \sum_{n=1}^{N} \frac{1}{n!} \underbrace{\int_0^t \cdots \int_0^t}_{n\ times} K_n(t,s_1,...,s_n)x(s_1)...x(s_n)ds_n \cdots ds_1 \qquad (1.2)$$

play a role analogous to ordinary polynomials. The classical theorem of Weierstrass, that continuous functions can be approximated by polynomials, has been generalized to nonlinear causal operators over functional spaces to the effect that general classes of nonlinear causal operators can be approximated by polynomial Volterra operators.
In this case the multiple nonlinear Volterra equation is a natural way to model a general nonlinear system with feedback. Referring to Figure 1 below, suppose that block $V$ of the system is described by equation (1.2), so that the output process $y(t)$ related to the input process

$$x(t) \text{ by } y(t) = V(x)(t), \tag{1.3}$$

where the operator $V$ is defined in (1.2). $u(t)$ is a control process, a known function. The feedback block $G$ is defined by the function $g(t, y(t))$, so that we have for the input process

$$x(t) = u(t) - g(t, y(t)). \tag{1.4}$$

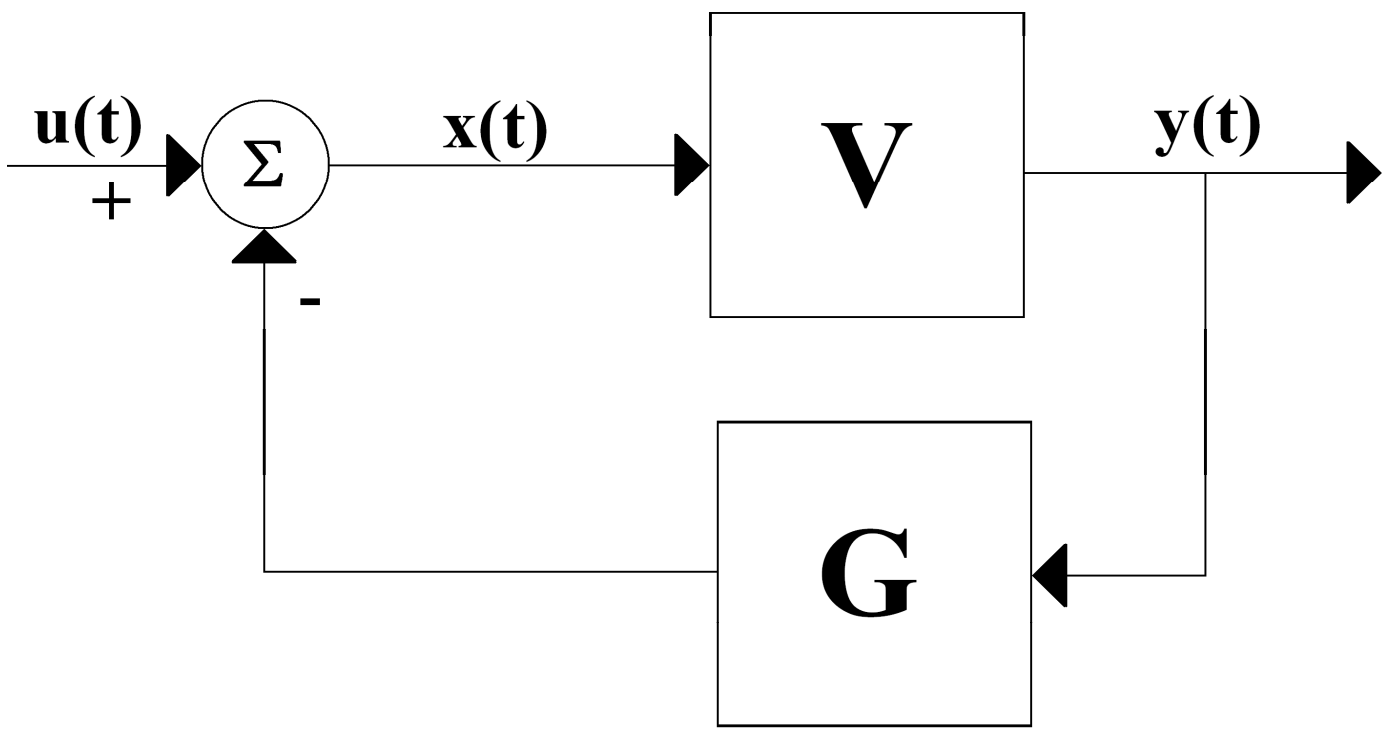

Fig.1. A General Nonlinear System with Feedback

Using equation (1.4) in (1.3), and taking into account (1.2), we obtain for $y(t)$:

$$y(t) = x_0(t) + \sum_{n=1}^{N} \frac{1}{n!} \underbrace{\int_0^t \cdots \int_0^t}_{n\ times} f_n(t, s_1, ..., s_n, y(s_1), ..., y(s_n)) ds_n \cdots ds_1, \tag{1.5}$$

where we have set

$$f_n(t, s_1, ..., s_n, y(s_1), ..., y(s_n)) := K_n(t, s_1, ..., s_n)(u(s_1) - g(s_1, y(s_1)))...(u(s_n) - g(s_n, y(s_n))).$$

Note that we just obtained in (1.5) one of the two equations that we are going to examine. The other is the infinite dimensional version below, which is an extension of (1.5):

$$y(t) = x_0(t) + \sum_{n=1}^{\infty} \frac{1}{n!} \underbrace{\int_0^t \cdots \int_0^t}_{n\ times} f_n(t, s_1, ..., s_n, y(s_1), ..., y(s_n)) ds_n \cdots ds_1. \tag{1.6}$$

The inclusion of the factor $\frac{1}{n!}$ in (1.6) is aimed at simplifying the conditions for convergence of the infinite series, and for consistency we have also included the same factor in (1.2).

Actually, our formulation covers not only the standard form of Volterra series, but also what we shall call generalized Volterra series; these are more general versions of what has been termed in [BCD] "Volterra-like series". Those "Volterra-like series" include terms of the form (in notation that is equivalent to the notation of [BCD])

$$\int_0^t \int_0^t \cdots \int_0^t K_M(t - s_1, t - s_2, ..., t - s_k)(x(s_1))^{m_1}(x(s_2))^{m_2} \cdots (x(s_k))^{m_k} ds_k \cdots ds_2 ds_1, \tag{1.7}$$

where $M := m_1 + m_2 + \cdots + m_k$. We note that [BCD] considers only time-invariant systems, and labels the kernels by the total degree of the corresponding monomial in $x$.

For our generalized Volterra series, we use standard multi-index notation. For a multi-index $\mathbf{m} := (m_1, m_2, ..., m_n) \in \mathbb{N}^n$, we set $|\mathbf{m}|_1 := \sum_{i=1}^{n} m_i$ and $\mathbf{m}! := \prod_{i=1}^{n} (m_i!)$. Then a generalized Volterra series has the form

$$\sum_{n=1}^{\infty} \sum_{\mathbf{m} \in \mathbb{N}^n : |\mathbf{m}|_1 \geq 1} \frac{1}{(n!)(\mathbf{m}!)} \int_0^t \cdots \int_0^t K_{n,\mathbf{m}}(t, s_1, ..., s_n) \left( \prod_{i=1}^{n} (y(s_i))^{m_i} \, ds_i \right) \tag{1.8}$$

Under the condition of uniform boundedness of the kernels $K_{n,\mathbf{m}}(t, s_1, ..., s_n)$, each series

$$\sum_{\mathbf{m} \in \mathbb{N}^n : |\mathbf{m}|_1 \geq 1} \frac{1}{(\mathbf{m}!)} K_{n,\mathbf{m}}(t, s_1, ..., s_n) \left( \prod_{i=1}^{n} (y(s_i))^{m_i} \right)$$

defines a real analytic function (analytic in the variables $y(s_i)$, $1 \leq i \leq n$), say $f_n(t, s_1, ..., s_n, y(s_1), ..., y(s_n))$, and, conversely, every such real analytic function $f_n$ , if its first derivatives with respect to each $y(s_i)$ do not vanish at $(0,..,0) \in \mathbb{R}^n$ can be, by definition, represented as a power series in the variables $y(s_i)$, $1 \leq i \leq n$, with the exponents $m_i$ in each monomial $\left( \prod_{i=1}^{n} (y(s_i))^{m_i} \right)$ being natural numbers (i.e. no exponent is zero). Of course, our problems, in this paper, are more general, since we employ conditions weaker than analyticity for the functions $f_n(t, s_1, ..., s_n, y(s_1), ..., y(s_n))$.

Furthermore, the models considered in this paper also include, as particular cases, equations of the form

$$x(t) = F(t, z_1, z_2, ...)$$

where each $z_i$ is a multiple Volterra-type integral of the form

$$\int_0^t \int_0^t \cdots \int_0^t g_k(t, s_1, s_2, ..., s_k, x(s_1), x(s_2), ..., x(s_k)) \, ds_k \cdots ds_2 \, ds_1$$

and $F$ is real analytic in the variables $z_1, z_2, ...$ .

Other applications come from modeling memory effects in physical, biological, economic, and social systems. In many instances, models with single Volterra equations are used only because there is reason to include memory effects; however, memory effects can result not only from single Volterra integral operators, but also from multiple Volterra integral operators. In addition, certain constitutive laws, for example in viscoelasticity, have been traditionally modelled via, among other possibilities, multiple Volterra integral operators.

We are interested in the numerical solution of (1.5) and (1.6). In recent years, there has been increasing interest in the numerical solution of many types of integral equations. However, the numerical solution of multiple Volterra equations, arising in the context of Volterra series or truncated Volterra series, has not been previously addressed in the literature. Ordinary Volterra equations have been studied, from the point of view of obtaining numerical solutions, in [H, L] and in other works. A problem of numerically solving a differential equation with an abstract memory term, that satisfies certain axiomatic conditions, has been studied in [CT].

The results of the present paper are not implied by any previous works. Furthermore, the proofs contained in the present paper include certain novel analytical and combinatorial arguments that have no counterpart in any previous work.

Our work in the present paper impinges on two general topics: the numerical solution of integral equations, and the representation of nonlinear in terms of (ordinary or generalized) Volterra series. The general topic of numerical solution of integral equations (not the type of equations treated in our paper, which constitute a novel type of problem formulated and studied here for the first time) has been, and continues to be, a topic of intensive research. As a sample of research activity in this area, we mention the papers [Ba, KS, RMA, TF]. The topic of applying Volterra series for modeling and identifying nonlinear dynamical systems is also a topic of vigorous contemporary research activity. Volterra series have found applications to a wide variety of real-world systems, including mechanical, electrical, biological, and economic and social systems. As instances of recent activity in this area, we cite the works [B, Ch1, Ch2]. An important topic in this area is the problem of identification of the input to a nonlinear system, modeled as a Volterra series, from observations of the output. This problem is too complex to briefly describe in this introduction, and we have devoted a separate section, section 3 of the present paper, to this problem and to its connection with the problems we have studied in the present paper.

A note about terminology: we use the term "multiple Volterra integral equations" to describe integral equations with multiple integrals in which the order of multiplicity of the integral exceeds the dimension of the independent variable. Integral equations in which the independent variable itself is multi-dimensional (as, for instance, in some applications of integral equations to multi-dimensional boundary value problems), but the order of the integrals is the same as the dimension of the independent variable (e.g. double integrals for two-dimensional boundary value problems, etc.) are not "multiple integral equations" in the context of our present paper. Similarly, we use the terms "finite dimensional" and "infinite dimensional" to refer to the multiplicity of the multiple integrals, not to the dimensions of the independent variables.

This paper is organized as follows: in section 2, we prove existence and uniqueness of solutions of the continuous-time problems (1.5) and (1.6) (these proofs are necessary before we start any discussion about numerical solutions); section 3 presents an application to the identification problem for nonlinear systems described by an input-output relationship in the form of a Volterra series; sections 4 and 5 deal with the numerical schemes for (1.5) and (1.6), respectively.

## 2. EXISTENCE AND UNIQUENESS OF SOLUTION BASED ON THE CONTRACTION MAPPING THEOREM

First, we want to state conditions guaranteeing that a solution $x(t)$ to equations (1.5) and (1.6) exists, and is unique in the space $X := C([0,T] \to \mathbb{R})$ of real-valued continuous functions on the interval $[0,T]$. $X$ is a complete normed linear space. Since $[0,T] \subset \mathbb{R}$ is compact, we may assume that a solution $x(t)$ of equation (1.5) is bounded, that is $|x(t)| \le R_x < \infty$, where $R_x$ depends on $x$. We use the notation

$\overline{B}^n(R) := \{(x_1, x_2, ..., x_n) : (\forall i = 1,2,...,n) : |x_i| \le R\}$, where $0 < R < \infty$.

We assume throughout this paper that, for $t \in [0,T]$, $x_0(t)$ is a known continuous differentiable function with bounded first derivative. We define

$$M_0 := \max_{t\in[0,T]} \left| \frac{dx_0(t)}{dt} \right|,$$

Also we set

$$C_0 := \max_{0 \le t \le T} |x_0(t)|.$$

Further we make the following assumptions:

(i) The kernels $f_1, f_2, ..., f_N$ are continuous in all their arguments, and, as a result, for each $n$ $f_n$ is bounded on compact sets. We set for each $n = 1,2..., N$:

$$C_n(R) := \max_{\substack{(t,s_1,...,s_n)\in\Delta_n([0,T]) \\ (\forall i): |x_i| \le R}} |f_n(t, s_1, ..., s_n; x_1, ..., x_n)| < \infty, \tag{2.1}$$

where $\Delta_n([0,T]) := \{(t, s_1, ..., s_n) : (\forall i = 1,..,n) : 0 \le s_i \le t \le T\}$. (2.2)

(ii) For each $n$, the kernel $f_n$ is Lipschitz in $(x_1, ..., x_n)$ on bounded sets, i.e. for each $n$ there exists $L_n(R) < \infty$, such that for each $(t, s_1, ..., s_n) \in \Delta_n([0,T])$ and for each $y, z \in \overline{B}^n(R)$ we have

$$|f_n(t, s_1, ..., s_n; y_1, ..., y_n) - f_n(t, s_1, ..., s_n; z_1, ..., z_n)| \le L_n(R) \sum_{i=1}^{n} |y_i - z_i|. \tag{2.3}$$

For every $Q > 0$, we define $\hat{L}_n(Q) := \sup\{L_n(R) : R \ge Q\}$.

(iii) We assume that

$$\overline{L}_n := \inf\{\hat{L}_n(Q) : Q > 0\} = \limsup_{R\to\infty} L_n(R) < \infty. \tag{2.4}$$

Theorem 2.1

*Under the stated conditions* (i), (ii), and (iii) *the N-dimensional Volterra equation* (1.5) *has a unique solution in* $X$ .

<u>Proof</u>

The proof of Theorem 2.1 will be similar to the one given in Corduneanu [1] for the case of one-dimensional Volterra equations. We introduce an operator $S^N$, acting on $X$, defined by

$$(S^N x)(t) = x_0(t) + \sum_{n=1}^{N} \frac{1}{n!} \underbrace{\int_0^t \cdots \int_0^t}_{n\ times} f_n(t, s_1, ..., s_n, x(s_1), ..., x(s_n)) ds_n \cdots ds_1 . \tag{2.5}$$

Observe that the right hand side of (2.5) is a continuous function on the interval $[0, T]$, in other words, $S^N$ maps $X$ into itself.

Next we define the weighted norm $\|\cdot\|_\mu$, where $\mu > 0$, on $X$ as follows

$$\| x \|_\mu := \max_{0 \le t \le T} (e^{-\mu t} | x(t) |) . \tag{2.6}$$

Then for any two functions $y(t)$ and $z(t)$ in $X$ we check that

$$e^{-\mu t} | (S^N y)(t) - (S^N z)(t) | \le \| y - z \|_\mu \frac{1 - e^{-\mu T}}{\mu} \sum_{n=1}^{N} \frac{\overline{L}_n}{(n-1)!} T^{n-1} . \tag{2.7}$$

Observe that the right hand side of (2.7) does not depend on variable $t$, hence

$$\| S^N y - S^N z \|_\mu \le \| y - z \|_\mu \frac{1 - e^{-\mu T}}{\mu} M^{(N)} , \tag{2.8}$$

where $M^{(N)} := \sum_{n=1}^{N} \frac{\overline{L}_n T^{n-1}}{(n-1)!}$ is a constant. Choosing $\mu$ sufficiently large, we make $q := \frac{1 - e^{-\mu T}}{\mu} M^{(N)} < 1$. Hence for the chosen value of $\mu$, we have

$$\| S^N y - S^N z \|_\mu \le q \| y - z \|_\mu , \text{ with } 0 < q < 1 . \tag{2.9}$$

Thus $S^N$ is a contraction on *X*, and therefore it has a unique fixed point $\overline{x}(t)$ in *X*; that fixed point is the unique solution of (1.6) in $X$ :

$$S^N(\overline{x}(t)) = (\overline{x}(t)) . \tag{2.10}$$

□

Note that the Contraction Mapping Theorem also gives the method of finding the solution. It is found by successive approximations: We choose an arbitrary $x^{(0)}(t) \in X$ , then

$$(\forall i): x^{(i)}(t) := S^N x^{(i-1)}(t), \text{ and } \bar{x}(t) = \lim_{i\to\infty} x^{(i)}(t).$$

For the above result to hold for the infinite-dimensional equation (1.6), we make two extra assumptions:

(iv) the following series is convergent for each $R < \infty$:

$$C_0 + \sum_{n=1}^{\infty} \frac{1}{n!} T^n C_n(R) < \infty, \tag{2.11}$$

(v) the following series with Lipschitz constants $\{\bar{L}_n\}_{n=1}^{\infty}$ in (2.4) converges:

$$M := \sum_{n=1}^{\infty} \frac{\bar{L}_n T^{n-1}}{(n-1)!} < \infty. \tag{2.12}$$

We have:

Theorem 2.2

*Under the stated conditions* (i), (ii), (iii), (iv), *and* (v), *the infinite-dimensional Volterra equation* (1.6) *has a unique solution in* $X$ .

Proof

Here we need to make only a few changes to the proof of Theorem 2.1. Now we consider the operator $S^{\infty}$ , acting on $X$ , and defined as follows:

$$(S^{\infty}x)(t) = x_0(t) + \sum_{n=1}^{\infty} \frac{1}{n!} \underbrace{\int_0^t \ldots \int_0^t}_{n \text{ times}} f_n(t, s_1, \ldots, s_n, x(s_1), \ldots, x(s_n)) ds_n \cdots ds_1. \tag{2.13}$$

Applying condition (iv), and Weierstrass M-test, we see that the series in the right hand side of (2.13) converges uniformly on $[0,T]$ to a continuous function, in other words, the operator $S^{\infty}$ maps the space $X$ into itself.

For any two functions $y(t)$ and $z(t)$ in $X$, we obtain

$$e^{-\mu t} |(S^{\infty}y)(t) - (S^{\infty}z)(t)| \le \| y - z \|_{\mu} \frac{1 - e^{-\mu T}}{\mu} M. \tag{2.14}$$

Since the right hand side of (2.14) does not depend on variable $t$, we have $\| S^{\infty}y - S^{\infty}z \|_{\mu} \leq q \| y - z \|_{\mu}$, where we select a value of $q := \frac{1 - e^{-\mu T}}{\mu} M$ strictly between 0 and 1 by choosing $\mu$ sufficiently large, so that $S^{\infty}$ is again a contraction on *X*, so it has a unique fixed point $\bar{x}(t)$ in *X*. This fixed point is a unique solution of (1.6) in $X$. □

## 3. APPLICATION TO MULTI-LINEAR VOLTERRA INTEGRAL EQUATIONS OF THE FIRST KIND

Multi-linear Volterra integral equations of the first kind arise naturally in the problem of input identification from observation of the output in general nonlinear systems modeled via input-ouput relationships in the formmof infinite-order Volterra series. Such problems have been extensively studied by Apartsin (Apartsyn) and Sidorov, see, e.g., [A1, A2, S1, S2].
In this section, we present a new method based on reduction (under suitable assumptions) of multi-linear Volterra integral equations of the first kind to multiple Volterra integral equations of infinite order of the second kind, of the type formulated in sections 1 and 2 of our previous paper. It is expected that the methods of numerical solution of infinite-order multiple Volterra integral equations of the second kind, as reported and analyzed in the following sections of our paper, will find practical applications to, among other things, the problem of identification of nonlinear systems modeled by multi-linear Volterra integral equations, thus adding new tools to the ones already studied by Apartsin and Sidorov and their collaborators.

In the theory of linear Volterra integral equations of the first kind, namely equations of the form

$$\int_0^t K(t,s)\,x(s)\,ds = f(t)\,, \tag{3.1}$$

there is a well known method for obtaining a solution, under suitable assumptions, by differentiating both sides of (3.1) with respect to *t*. This method is due to Volterra himself. The result of differentiating (3.1) with respect to *t* is

$$K(t,t)x(t) + \int_0^t K_t(t,s)\,x(s)\,ds = f_t(t) \tag{3.2}$$

and, assuming that $K(t, t)$ never vanishes for *t* in some interval $[0,T]$, we obtain from (3.2) the second-kind Volterra integral equation

$$x(t) = \frac{f_t(t)}{K(t,t)} - \int_0^t \frac{K_t(t,s)}{K(t,t)}\,x(s)\,ds\,. \tag{3.3}$$

We are interested in extending the same idea to multi-linear Volterra integral equations of the first kind. Such equations arise naturally in identification problems for systems described by Volterra series, and their analytical properties and numerical solution have been studied primarily by Apartsin, Sidorov, and their collaborators. A multi-linear Volterra integral equation of the first kind has the form

$$\sum_{n=1}^{\infty} \frac{1}{n!} \int_{s_1=0}^{t} \int_{s_2=0}^{t} \cdots \int_{s_n=0}^{t} K_n(t,s_1,s_2,\ldots,s_n)\,x(s_1)\,x(s_2)\cdots x(s_n)\,ds_1\,ds_2\cdots ds_n = f(t) \tag{3.4}$$

Eq. (3.4) arises in the problem of identification of the input to a nonlinear system, if the output is observed and the system is described by a Volterra series with known kernels.

One of the difficulties of numerically solving (3.4) stems from the fact that, by using a direct discretization of (3.4), we get a system of polynomial equations in the unknown values of the approximation to $x(t)$ at the discretization nodes. This is to be contrasted with the case of second-order multi-linear Volterra integral equations, where discretization leads to recursive formulas for the wanted approximations, so that the solution of the discretized equations presents no theoretical difficulty. For multi-linear equations of low order (for example, bilinear, trilinear, etc.), it may be possible to get some analytical results for existence and uniqueness of solutions of the system of polynomial equations, as well as guaranteed algorithms for the solution, but for general finite-order multi-linear equations, and for infinite-order multi-linear equations, that does not seem to be possible in general. Therefore, we look for an extension of Volterra's method to reduce the first-kind multi-linear equation to a second-kind multi-linear equation. By differentiating both sides of (3.4) with respect to *t*, and assuming without loss of generality that each $K_n(t,s_1,s_2,...,s_n)$ is symmetric with respect to $s_1,s_2,...,s_n$, we find

$$\left\{K_1(t,t)+\sum_{n=2}^{\infty}\frac{1}{(n-1)!}\int_{s_1=0}^{t}\int_{s_2=0}^{t}\cdots\int_{s_{n-1}=0}^{t}K_n(t,s_1,\ldots,s_{n-1},t)x(s_1)x(s_2)\cdots x(s_{n-1})\,ds_{n-1}\cdots ds_2\,ds_1\right\}x(t)+$$
$$+\sum_{n=1}^{\infty}\frac{1}{n!}\int_{s_1=0}^{t}\int_{s_2=0}^{t}\cdots\int_{s_{n-1}=0}^{t}K_{n,t}(t,s_1,\ldots,s_{n-1},s_n)x(s_1)x(s_2)\cdots x(s_n)\,ds_n\cdots ds_2\,ds_1=f_t(t)\ . \tag{3.5}$$

Assuming that the expression in curly brackets in (3.5) never vanishes, we can solve (3.5) for $x(t)$. For simplicity, we set

$$J(t):=K_1(t,t)+\sum_{n=2}^{\infty}\frac{1}{(n-1)!}\int_{s_1=0}^{t}\int_{s_2=0}^{t}\cdots\int_{s_{n-1}=0}^{t}K_n(t,s_1,\ldots,s_{n-1},t)x(s_1)x(s_2)\cdots x(s_{n-1})\,ds_{n-1}\cdots ds_2\,ds_1 \tag{3.6}$$

Then we have

$$x(t)=\frac{f_t(t)}{J(t)}-\frac{1}{J(t)}\sum_{n=1}^{\infty}\frac{1}{n!}\int_{s_1=0}^{t}\int_{s_2=0}^{t}\cdots\int_{s_{n-1}=0}^{t}K_{n,t}(t,s_1,\ldots,s_{n-1},s_n)x(s_1)x(s_2)\cdots x(s_n)\,ds_n\cdots ds_2\,ds_1\,. \tag{3.7}$$

We set

$$L_n(t,s_1,\cdots,s_n):=\frac{K_{n+1}(t,s_1,\cdots,s_n,t)}{K_1(t,t)}, \tag{3.8}$$

$$G(t):=1+\sum_{n=1}^{\infty}\frac{1}{n!}\int_{s_1=0}^{t}\int_{s_2=0}^{t}\cdots\int_{s_{n-1}=0}^{t}L_n(t,s_1,\ldots,s_{n-1},s_n)x(s_1)x(s_2)\cdots x(s_n)\,ds_n\cdots ds_2\,ds_1 \tag{3.9}$$

and we rewrite (3.7 ) as

$$x(t)=\frac{1}{G(t)}\frac{f_t(t)}{K_1(t,t)}-\sum_{n=1}^{\infty}\frac{1}{(n!)G(t)}\int_{s_1=0}^{t}\int_{s_2=0}^{t}\cdots\int_{s_{n-1}=0}^{t}\frac{K_{n,t}(t,s_1,\ldots,s_{n-1},s_n)}{K_1(t,t)}x(s_1)x(s_2)\cdots x(s_n)\,ds_n\cdots ds_2\,ds_1 \tag{3.10}$$

Further, assuming that

$$\left|\sum_{n=1}^{\infty}\frac{1}{n!}\int_{s_1=0}^{t}\int_{s_2=0}^{t}\cdots\int_{s_{n-1}=0}^{t}L_n(t,s_1,\ldots,s_{n-1},s_n)x(s_1)x(s_2)\cdots x(s_n)\,ds_n\cdots ds_2\,ds_1\right|<1\ , \tag{3.11}$$

we expand $\frac{1}{G(t)}$ as follows:

$$\frac{1}{G(t)}=1+\sum_{m=1}^{\infty}(-1)^m\left\{\sum_{n=1}^{\infty}\frac{1}{n!}\int_{s_1=0}^{t}\int_{s_2=0}^{t}\cdots\int_{s_{n-1}=0}^{t}L_n(t,s_1,\ldots,s_{n-1},s_n)x(s_1)x(s_2)\cdots x(s_n)\,ds_n\cdots ds_2\,ds_1\right\}^m . \tag{3.12}$$

Before proceeding with further calculations, we introduce some terminology and notation. If $N$ is a natural number, an <u>ordered partition</u> of $N$ is an ordered collection of natural numbers $(n_1,n_2,\ldots,n_k)$ such that $n_1+n_2+\cdots+n_k=N$. We denote the collection of all ordered partitions of $N$ by $\mathbf{P}_N$. A particular partition of $N$ will be denoted by $P$, and the number $k$ of natural numbers in the partition $P$ will be denoted by $\|P\|$.

With this notation and terminology, we have

$$\left\{\sum_{n=1}^{\infty}\frac{1}{n!}\int_{s_1=0}^{t}\int_{s_2=0}^{t}\cdots\int_{s_n=0}^{t}L_n(t,s_1,\ldots,s_{n-1},s_n)x(s_1)x(s_2)\cdots x(s_n)\,ds_n\cdots ds_2\,ds_1\right\}^m=$$

$$=\sum_{N=1}^{\infty}\sum_{\substack{(n_1,n_2,\ldots,n_{\|P\|})\equiv P\in\mathbf{P}_N\\ \|P\|=m}}\left(\prod_{i=1}^{m}\frac{1}{n_i!}\right)\int_{s_1=0}^{t}\int_{s_2=0}^{t}\cdots\int_{s_N=0}^{t}L_{n_1}(t,s_1,\cdots,s_{n_1})L_{n_2}(t,s_{n_1+1},\cdots,s_{n_1+n_2})\cdots \tag{3.13}$$

$$\cdot L_{n_{\|P\|}}(t,s_{n_1+n_2+\cdots+n_{\|P\|-1}},\cdots,s_N)\,x(s_1)x(s_2)\cdots x(s_N)\,ds_N\cdots ds_2\,ds_1\ .$$

Eq. (3.13) is proved as follows:

We start with the well known theorem about the product of two absolutely convergent series:

$$\left(\sum_{n=1}^{\infty} A_n\right)\left(\sum_{n=1}^{\infty} B_n\right)=\sum_{n=1}^{\infty}\sum_{k=1}^{n} A_k B_{n-k}\,. \tag{3.14}$$

By repeated application of (3.14) above, we find that the product of $m$ absolutely convergent series, say the series $\sum_{n=1}^{\infty} A_{(i)n}$ , $i = 1,2,...,m$ , is expressed as

$$\prod_{i=1}^{m}\left(\sum_{n=1}^{\infty} A_{(i)n}\right)=\sum_{n=1}^{\infty}\sum_{\substack{(k_1,k_2,...,k_m):\\ k_1+k_2+\cdots+k_m=n}}\prod_{i=1}^{m} A_{(i)k_i} \tag{3.15}$$

and therefore, for the the $m$-th power of an absolutely convergent series (i.e. a product of the type shown on the left-hand side of (3. ) above but with all series identical to each other), we have

$$\left(\sum_{n=1}^{\infty} A_n\right)^m=\sum_{n=1}^{\infty}\sum_{\substack{(k_1,k_2,...,k_m):\\ k_1+k_2+\cdots+k_m=n}}\prod_{i=1}^{m} A_{k_i}\,. \tag{3.16}$$

By taking each $A_n$ to be equal to

$$\frac{1}{n!}\int_{s_1=0}^{t}\int_{s_2=0}^{t}\cdots\int_{s_n=0}^{t} L_n(t,s_1,\ldots,s_{n-1},s_n)x(s_1)x(s_2)\cdots x(s_n)\,ds_n\cdots ds_2\,ds_1$$

we find

$$\left\{\sum_{n=1}^{\infty}\frac{1}{n!}\int_{s_1=0}^{t}\int_{s_2=0}^{t}\cdots\int_{s_n=0}^{t} L_n(t,s_1,\ldots,s_{n-1},s_n)x(s_1)x(s_2)\cdots x(s_n)\,ds_n\cdots ds_2\,ds_1\right\}^m =$$

$$=\sum_{N=1}^{\infty}\sum_{\substack{(n_1,n_2,\ldots,n_{\|P\|})\equiv P\in\mathbf{P}_N\\ \|P\|=m}}\prod_{i=1}^{m}\left(\frac{1}{n_i!}\int_{s_1=0}^{t}\int_{s_2=0}^{t}\cdots\int_{s_N=0}^{t} L_{n_i}(t,s_1,\cdots,s_{n_i})\;x(s_1)x(s_2)\cdots x(s_{n_i})\,ds_{n_i}\cdots ds_2\,ds_1\right) \tag{3.17}$$

and then, by re-naming the variables in the multiple integrals on the right-hand side of (3.17), so that

$$\prod_{i=1}^{m}\left(\frac{1}{n_i!}\int_{s_1=0}^{t}\int_{s_2=0}^{t}\cdots\int_{s_N=0}^{t} L_{n_i}(t,s_1,\cdots,s_{n_i})\, x(s_1)x(s_2)\cdots x(s_{n_i})\, ds_{n_i}\cdots ds_2\, ds_1\right)=$$

$$=\left(\prod_{i=1}^{m}\frac{1}{n_i!}\right)\int_{s_1=0}^{t}\int_{s_2=0}^{t}\cdots\int_{s_N=0}^{t} L_{n_1}(t,s_1,\cdots,s_{n_1})\, L_{n_2}(t,s_{n_1+1},\cdots,s_{n_1+n_2})\cdots \tag{3.18}$$

$$\cdot L_{n_{\|P\|}}(t,s_{n_1+n_2+\cdots+n_{\|P\|-1}},\cdots,s_N)\, x(s_1)x(s_2)\cdots x(s_N)\, ds_N\cdots ds_2\, ds_1\ ,$$

we obtain (3.13). This completes the proof of (3.13).

Consequently, we can write

$$\frac{1}{G(t)}=1+\sum_{n=1}^{\infty}\int_{s_1=0}^{t}\int_{s_2=0}^{t}\cdots\int_{s_n=0}^{t} M_n(t,s_1,\ldots,s_{n-1},s_n)x(s_1)x(s_2)\cdots x(s_n)\, ds_n\cdots ds_2\, ds_1; \tag{3.19}$$

where

$$M_n(t,s_1,\ldots,s_{n-1},s_n):=\sum_{m=1}^{\infty}(-1)^m\sum_{\substack{(n_1,n_2,\ldots,n_{\|P\|})\equiv P\in\mathbf{P}_n\\ \|P\|=m}}\left(\prod_{i=1}^{m}\frac{1}{n_i!}\right)L_{n_1}(t,s_1,\cdots,s_{n_1})\, L_{n_2}(t,s_{n_1+1},\cdots,s_{n_1+n_2})\cdots$$

$$\cdot L_{n_{\|P\|}}(t,s_{n_1+n_2+\cdots+n_{\|P\|-1}},\cdots,s_N)\ . \tag{3.20}$$

Furthermore, we have

$$\left(1+\sum_{n=1}^{\infty}\int_{s_1=0}^{t}\int_{s_2=0}^{t}\cdots\int_{s_n=0}^{t} M_n(t,s_1,\ldots,s_{n-1},s_n)x(s_1)x(s_2)\cdots x(s_n)\, ds_n\cdots ds_2\, ds_1\right)\cdot$$

$$\cdot\left(\frac{f_t(t)}{K_1(t,t)}-\sum_{n=1}^{\infty}\frac{1}{n!}\int_{s_1=0}^{t}\int_{s_2=0}^{t}\cdots\int_{s_n=0}^{t}\frac{K_{n,t}(t,s_1,\ldots,s_{n-1},s_n)}{K_1(t,t)}x(s_1)x(s_2)\cdots x(s_n)\, ds_n\cdots ds_2\, ds_1\right)=$$

$$=\frac{f_t(t)}{K_1(t,t)}+\sum_{n=1}^{\infty}\int_{s_1=0}^{t}\int_{s_2=0}^{t}\cdots\int_{s_n=0}^{t} Q_n(t,s_1,\ldots,s_{n-1},s_n)x(s_1)x(s_2)\cdots x(s_n)\, ds_n\cdots ds_2\, ds_1 \tag{3.21}$$

where

$$Q_n(t,s_1,\ldots,s_{n-1},s_n) := \frac{f_t(t)}{K_1(t,t)} M_n(t,s_1,\ldots,s_{n-1},s_n) - \frac{1}{n!}\frac{K_{n,t}(t,s_1,\ldots,s_n)}{K_1(t,t)} - \\ -\sum_{k=1}^{n} \frac{1}{k!}\frac{K_{k,t}(t,s_1,\ldots,s_k)}{K_1(t,t)} M_{n-k}(t,s_{k+1},\ldots,s_n) \quad . \tag{3.22}$$

In this way, we obtain the multiple Volterra equation of the second kind

$$x(t) = \frac{f_t(t)}{K_1(t,t)} + \sum_{n=1}^{\infty} \int_{s_1=0}^{t} \int_{s_2=0}^{t} \cdots \int_{s_n=0}^{t} Q_n(t,s_1,\ldots,s_{n-1},s_n) x(s_1) x(s_2) \cdots x(s_n)\, ds_n \cdots ds_2\, ds_1 \ . \tag{3.23}$$

This is precisely a multiple Volterra equation of the type studied in this paper.

4. FINITE DIFFERENCE SCHEME FOR THE N-DIMENSIONAL CASE

We take a partition of the interval $[0,T]$ with step $h$: let $t_i := ih,\ i = 0,1,...,M$ : $t_0 = 0, t_0 < t_1 < ... < t_M = T$ . We approximate the values $\{x(t_i)\}_{i=0}^{M}$ of the exact solution of equation (1.5)on the partition $\{t_i\}_{i=0}^{M}$, by the solution $\{x_i\}_{i=0}^{M}$ of the following finite-difference equation:

$$x_i = x_0(t_i) + \sum_{n=1}^{N} \frac{h^n}{n!} \sum_{I_n(i)} f_n(t_i; t_{j_1}, \dots, t_{j_n}; x_{j_1}, \dots, x_{j_n}), \tag{4.1}$$

where

$$I_n(i) := \{(j_1, j_2, \dots, j_n) : (\forall k = 1,2,\dots,n) : 0 \le j_k \le i-1\}. \tag{4.2}$$

We denote the error of this approximation at $t_i$ by

$$\varepsilon_i := |x(t_i) - x_i|. \tag{4.3}$$

<u>Theorem 4.1</u>
*We assume that*
(i*) the kernels* $f_n$ *in (1.5) are differentiable both with respect to* $\vec{s} := (s_1, \dots, s_n)$ *and with respect to* $t$*, moreover*

$$|\nabla_{\vec{s}} f_n(t_i, s_1, \dots, s_n, x_1, \dots, x_n)| \le M_n, \tag{4.4}$$

$$\left|\frac{\partial f_n}{\partial t}(t, s_1, \dots, s_n, x_1, \dots, x_n)\right| \le B_n; \tag{4.5}$$

(ii) *the Lipschitz condition holds for the kernels in (1.5) for each n:*

$$|f_n(t, s_1, \dots, s_n; z_1, \dots, z_n) - f_n(t, s_1, \dots, s_n; y_1, \dots, y_n)| \le L_n\{|z_1 - y_1| + \dots + |z_n - y_n|\}; \tag{4.6}$$

(iii) *the kernels are bounded:*

$$|f_n(t, s_1, \dots, s_n, z_1, \dots, z_n)| \le C_n. \tag{4.7}$$

*Then for all* $i$ *the absolute value of error* $|\varepsilon_i|$ *of the numerical approximation (4.1)-(4.2) of the solution of equation (1.5) is of the order of* $h$*, uniformly in N.*

<u>Proof</u>
With the notation in place we can represent $x(t_i)$ in the following form:

$$x(t_i) = x_0(t_i) + \sum_{n=1}^{N} \frac{1}{n!} \sum_{I_n(i)} \underbrace{\int_{t_{j_1}}^{t_{j_1+1}} \cdots \int_{t_{j_n}}^{t_{j_n+1}}}_{n \text{ times}} f_n(t_i, s_1, \dots, s_n; x(s_1), \dots, x(s_n)) ds_n \cdots ds_1. \tag{4.8}$$

Then for the absolute value of the error we have

$$| \varepsilon_i | \le \sum_{n=1}^{N} \frac{1}{n!} \sum_{I_n(i)} \int_{t_{j_1}}^{t_{j_1+1}} \cdots \int_{t_{j_n}}^{t_{j_n+1}} | f_n(t_i, s_1,..., s_n; x(s_1),..., x(s_n)) -$$
$$- f_n(t_i, t_{j_1},..., t_{j_n}; x_{j_1},..., x_{j_n}) | ds_n \cdots ds_1 \ . \quad (4.9)$$

We can represent the expression inside the absolute value bars in (4.9) as follows:

$$f_n(t_i, s_1,..., s_n; x(s_1),..., x(s_n)) - f_n(t_i, t_{j_1},..., t_{j_n}; x_{j_1},..., x_{j_n}) =$$
$$= (f_n(t_i, s_1,..., s_n; x(s_1),..., x(s_n)) - f_n(t_i, t_{j_1},..., t_{j_n}; x(s_1),..., x(s_n))) +$$
$$+ (f_n(t_i, t_{j_1},..., t_{j_n}; x(s_1),..., x(s_n)) - f_n(t_i, t_{j_1},..., t_{j_n}, x(t_{j_1}),..., x(t_{j_n}))) +$$
$$+ (f_n(t_i, t_{j_1},..., t_{j_n}; x(t_{j_1}),..., x(t_{j_n})) - f_n(t_i, t_{j_1},..., t_{j_n}; x_{j_1},..., x_{j_n}) \ . \quad (4.10)$$

Using (4.10) in (4.9), we obtain an upper bound for the absolute value of error:

$$| \varepsilon_i | \le \sum_{n=1}^{N} \frac{1}{n!} \sum_{I_n(i)} \int_{t_{j_1}}^{t_{j_1+1}} \cdots \int_{t_{j_n}}^{t_{j_n+1}} | f_n(t_i, s_1,..., s_n; x(s_1),..., x(s_n)) -$$
$$- f_n(t_i, t_{j_1},..., t_{j_n}, x(s_1),..., x(s_n)) | ds_n ... ds_1 +$$
$$+ \sum_{n=1}^{N} \frac{1}{n!} \sum_{I_n(i)} \int_{t_{j_1}}^{t_{j_1+1}} \cdots \int_{t_{j_n}}^{t_{j_n+1}} | f_n(t_i, t_{j_1},..., t_{j_n}; x(s_1),..., x(s_n)) -$$
$$- f_n(t_i, t_{j_1},..., t_{j_n}; x(t_{j_1}),..., x(t_{j_n})) | ds_n ... ds_1 +$$
$$+ \sum_{n=1}^{N} \frac{1}{n!} \sum_{I_n(i)} \int_{t_{j_1}}^{t_{j_1+1}} \cdots \int_{t_{j_n}}^{t_{j_n+1}} | f_n(t_i, t_{j_1},..., t_{j_n}; x(t_{j_1}),..., x(t_{j_n})) -$$
$$- f_n(t_i, t_{j_1},..., t_{j_n}; x_{j_1},..., x_{j_n}) | ds_n ... ds_1 \ . \quad (4.11)$$

For each $k = 1,..., n$ in (4.11) we have $t_{j_k} \le s_k \le t_{j_k+1} = t_{j_k} + h$, hence

$$\| (s_1,..., s_n) - (t_{j_1},..., t_{j_n}) \|_{\mathbb{R}^n} \le h\sqrt{n} \ . \quad (4.12)$$

By the Mean Value Theorem, for some $\vartheta \in (0,1)$, we have

$$| f_n(t_i, s_1,..., s_n; z_1,..., z_n) - f_n(t_i, t_{j_1},..., t_{j_n}; z_1,..., z_n) | =$$
$$= | \nabla_{\vec{s}} f_n(t_i, (1-\vartheta)s_1 + \vartheta t_{j_1},..., (1-\vartheta)s_n + \vartheta t_{j_n}; z_1,..., z_n) \cdot \vec{u} \, | h\sqrt{n} \ , \quad (4.13)$$

where $\vec{u}$ is the unit vector in $\mathbb{R}^n$ directed from $(s_1,..., s_n)$ to $(t_{j_1},..., t_{j_n})$.
By assumption (i), we have

$$| f_n(t_i, s_1,..., s_n; z_1,..., z_n) - f_n(t_i, t_{j_1},..., t_{j_n}; z_1,..., z_n) | \le M_n h\sqrt{n} \ , \quad (4.14)$$

and for the first term in the right hand side in (4.11) we have

$$\sum_{n=1}^{N}\frac{1}{n!}\sum_{I_n(i)}\int_{t_{j_1}}^{t_{j_1+1}}\cdots\int_{t_{j_n}}^{t_{j_n+1}}|f_n(t_i,s_1,...,s_n;x(s_1),...,x(s_n))-$$

$$-f_n(t_i,t_{j_1},...,t_{j_n};x(s_1),...,x(s_n))|\,ds_n...ds_1\le$$

$$\le\sum_{n=1}^{N}\frac{1}{n!}\sum_{I_n(i)}\int_{t_{j_1}}^{t_{j_1+1}}\cdots\int_{t_{j_n}}^{t_{j_n+1}}M_n h\sqrt{n}\,ds_n...ds_1\le hE^{(N)}=O(h) \tag{4.15}$$

where we have used the notation

$$E^{(N)}:=\sum_{n=1}^{N}\frac{M_n\sqrt{n}T^n}{n!}. \tag{4.16}$$

Applying (ii), we get an upper bound for the third term in the right hand side of (4.11):

$$\sum_{n=1}^{N}\frac{1}{n!}\sum_{I_n(i)}\int_{t_{j_1}}^{t_{j_1+1}}\cdots\int_{t_{j_n}}^{t_{j_n+1}}|f_n(t_i,t_{j_1},...,t_{j_n},x(t_{j_1}),...,x(t_{j_n}))-f_n(t_i,t_{j_1},...,t_{j_n},x_{j_1},...,x_{j_n})|\,ds_n...ds_1\le$$

$$\le\sum_{n=1}^{N}\frac{L_n h^n}{n!}\sum_{I_n(i)}\sum_{k=1}^{n}|\varepsilon_{j_k}| \tag{4.17}$$

Suppose $t_2>t_1$. Using conditions (i) and (iii), we obtain

$$|x(t_2)-x(t_1)|\le M_0(t_2-t_1)+\sum_{n=1}^{N}\frac{1}{n!}\left[B_n(t_2-t_1)t_1^n+C_n(t_2^n-t_1^n)\right] \tag{4.18}$$

Now, using condition (ii), we have for the second term in the right hand side of (4.15):

$$\sum_{n=1}^{N}\frac{1}{n!}\sum_{I_n(i)}\int_{t_{j_1}}^{t_{j_1+1}}\cdots\int_{t_{j_n}}^{t_{j_n+1}}|f_n(t_i,t_{j_1},...,t_{j_n},x(s_1),...,x(s_n))-f_n(t_i,t_{j_1},...,t_{j_n},x(t_{j_1}),...,x(t_{j_n}))|\,ds_n...ds_1\le$$

$$\le\sum_{n=1}^{N}\frac{L_n h^{n-1}}{n!}\sum_{I_n(i)}\sum_{k=1}^{n}\int_{t_{j_k}}^{t_{j_k+1}}\left|x(s_k)-x(t_{j_k})\right|ds_k. \tag{4.19}$$

For the integral in the right hand side of (4.19) we have, using (4.18)

$$\int_{t_{j_k}}^{t_{j_k+1}}\left|x(s_k)-x(t_{j_k})\right|ds_k\le\frac{M_0h^2}{2}+\sum_{m=1}^{N}\frac{1}{m!}\left[B_m t_{j_k}^m\frac{h^2}{2}+C_m\frac{(t_{j_k}+h)^{m+1}-t_{j_k}^{m+1}-(m+1)t_{j_k}^m h}{m+1}\right]$$

(4.20)

Using the Binomial Theorem we rewrite right hand side of (4.20) as follows:

$$\frac{M_0 h^2}{2} + \sum_{m=1}^{N} \left[ \frac{B_m t_{j_k}^m h^2}{2m!} + \frac{C_m}{(m+1)!} \left\{ \sum_{l=2}^{m+1} \binom{m+1}{l} t_{j_k}^{m+1-l} h^l \right\} \right]. \tag{4.21}$$

In the inside summation in (4.21) terms, corresponding to $l > 2$, are of higher order of error; for $l = 2$ we have:

$$\frac{C_m}{(m+1)!} \binom{m+1}{2} t_{j_k}^{m-1} h^2 = \frac{C_m}{2(m-1)!} t_{j_k}^{m-1} h^2 . \tag{4.22}$$

Using (4.21) and (4.22), we obtain for (4.20):

$$\int_{t_{jk}}^{t_{jk+1}} \left| x(s_k) - x(t_{j_k}) \right| ds_k \le \left\{ \frac{M_0}{2} + \sum_{m=1}^{N} \left[ \frac{B_m t_{j_k}^m}{2m!} + \frac{C_m t_{j_k}^{m-1}}{2(m-1)!} \right] \right\} h^2 . \tag{4.23}$$

Making use of (4.23) in the right hand side of (4.19), we obtain

$$\sum_{n=1}^{N} \frac{L_n h^{n-1}}{n!} \sum_{I_n(i)} \sum_{k=1}^{n} \int_{t_{jk}}^{t_{jk+1}} \left| x(s_k) - x(t_{j_k}) \right| ds_k \le$$

$$\le \sum_{n=1}^{N} \left( \frac{L_n h^{n+1}}{2n!} \sum_{I_n(i)} \frac{n M_0}{2} \right) + \sum_{n=1}^{N} \left\{ \frac{L_n h^{n+1}}{2n!} \sum_{I_n(i)} \sum_{m=1}^{N} \sum_{k=1}^{n} \left( \frac{B_m t_i^m}{m!} + \frac{C_m t_i^{m-1}}{(m-1)!} \right) \right\}. \tag{4.24}$$

For the first term in the right hand side of (4.24) we have

$$\sum_{n=1}^{N} \left( \frac{L_n h^{n+1}}{2n!} \sum_{I_n(i)} \frac{n M_0}{2} \right) \le \frac{M_0 D^{(N)} T}{2} h \tag{4.25}$$

where we denoted

$$D^{(N)} := \sum_{n=1}^{N} \frac{L_n T^{n-1}}{(n-1)!} . \tag{4.26}$$

For the second term in the right hand side of (4.24) we have

$$\sum_{n=1}^{N} \frac{L_n h^{n+1}}{2n!} \sum_{I_n(i)} \sum_{m=1}^{N} \sum_{k=1}^{n} \left( \frac{B_m t_i^m}{m!} + \frac{C_m t_i^{m-1}}{(m-1)!} \right) \le \frac{(B^{(N)} + C^{(N)}) D^{(N)} T}{2} h \tag{4.27}$$

where we denoted

$$B^{(N)} := \sum_{m=1}^{N} \frac{B_m T^m}{m!}, \tag{4.28}$$

and

$$C^{(N)} := \sum_{m=1}^{N} \frac{C_m T^{m-1}}{(m-1)!}. \tag{4.29}$$

Using (4.27) and (4.25) in (4.24), we obtain an upper bound for the second term (4.11):

$$\sum_{n=1}^{N} \frac{1}{n!} \sum_{I_n(i)} \int_{t_{j_1}}^{t_{j_1+1}} \cdots \int_{t_{j_n}}^{t_{j_n+1}} | f_n(t_i, t_{j_1}, ..., t_{j_n}, x(s_1), ..., x(s_n)) - f_n(t_i, t_{j_1}, ..., t_{j_n}, x(t_{j_1}), ..., x(t_{j_n})) | \, ds_n ... ds_1 \le$$

$$\le \frac{(M_0 + B^{(N)} + C^{(N)}) D^{(N)} T}{2} h. \tag{4.30}$$

At this point, we need to prove the following lemma:

<u>Lemma 4.1</u>
*The following equality holds:*

$$\sum_{(j) \in I_n(i)} \sum_{k=1}^{n} \delta_{j_k} = n i^{n-1} \sum_{l=0}^{i-1} \delta_l \ . \tag{4.31}$$

<u>Proof of lemma 4.1</u>
Using multi-index notation $(l) := (l_0, l_1, ..., l_{i-1})$, we define

$$L_n^{(i)} := \{(l) : (\forall k = 0,1,...,i-1) : l_k \in \{0,1,...,n\} \text{ and } \sum_{k=0}^{i-1} l_k = n \ \}, \tag{4.32}$$

and

$$S_{(l)}^{(i)}(n) := \{(j) \in I_n(i) : (\forall m = 0,1,...,i-1) : \#(k : j_k = m) = l_k\}. \tag{4.33}$$

Then we can partition $I_n(i)$ as follows:

$$I_n(i) = \mathop{\dot{\bigcup}}_{(l) \in L_n^{(i)}} S_{(l)}^{(i)}(n). \tag{4.34}$$

Note that

$$|S_{(l)}^{(i)}(n)| = \binom{n}{l_0\, l_1 \ldots l_{i-1}}. \tag{4.35}$$

With the new notation in place we have:

$$\sum_{(j)\in I_n(i)} \sum_{k=1}^{n} \delta_{j_k} = \sum_{(l)\in L_n^{(i)}} \sum_{(j)\in S_{(l)}^{(i)}(n)} \sum_{k=1}^{n} \delta_{j_k} = \sum_{(l)\in L_n^{(i)}} \sum_{(j)\in S_{(l)}^{(i)}(n)} \left( \sum_{j=0}^{i-1} l_j \delta_j \right) =$$
$$= \sum_{(l)\in L_n^{(i)}} |S_{(l)}^{(i)}(n)| \left( \sum_{j=0}^{i-1} l_j \delta_j \right) = \sum_{(l)\in L_n^{(i)}} \binom{n}{l_0\, l_1 \ldots l_{i-1}} \left( \sum_{j=0}^{i-1} l_j \delta_j \right) =$$
$$= \sum_{j=0}^{i-1} \sum_{(l)\in L_n^{(i)}} \binom{n}{l_0\, l_1 \ldots l_{i-1}} l_j \delta_j = n \sum_{j=0}^{i-1} \sum_{(l)'_j \in L_{n-1}^{(i)}} \binom{n-1}{l_0\, l_1 \ldots l'_j \ldots l_{i-1}} \delta_j , \tag{4.36}$$

where $(l)'_j := (l_0, l_1, \ldots, l_{j-1}, l'_j, l_{j+1}, \ldots, l_{i-1})$, and $l'_j := l_j - 1$.

Now, using the Multinomial Theorem, we rewrite the right hand side of (4.36):

$$n \sum_{j=0}^{i-1} \sum_{(l)'_j \in L_{n-1}^{(i)}} \binom{n-1}{l_0\, l_1 \ldots l'_j \ldots l_{i-1}} \delta_j = n \sum_{j=0}^{i-1} \underbrace{(1+1+\ldots+1)}_{i \text{ terms}}^{n-1} \delta_j = n \sum_{j=0}^{i-1} i^{n-1} \delta_j = n i^{n-1} \sum_{j=0}^{i-1} \delta_j ,$$

as required. □

We now continue the proof of theorem 4.1.

Using (4.31), we obtain an upper bound for the right hand side of (4.17):

$$\sum_{n=1}^{N} \frac{L_n h^n}{n!} \sum_{I_n(i)} \sum_{k=1}^{n} |\varepsilon_{j_k}| \le D^{(N)} h \sum_{l=0}^{i-1} |\varepsilon_l| . \tag{4.37}$$

Combining (4.15), (4.17), (4.30), and (4.37), we obtain an upper bound for the absolute value of the error $\varepsilon_i$ in (4.11):

$$|\varepsilon_i| \le F^{(N)} h + D^{(N)} h \sum_{j=0}^{i-1} |\varepsilon_j| , \tag{4.38}$$

where we denote

$$F^{(N)} := \frac{D^{(N)} T (M_0 + B^{(N)} + C^{(N)})}{2} + E^{(N)} . \tag{4.39}$$

We note that inequality (4.38) is the discrete Gronwall inequality stated in the lemma below, with $y_i = |\varepsilon_i|$, $a = F^{(N)}h$, $b = D^{(N)}h$, and $y_0 = |\varepsilon_0| = 0$.
We shall need the following

<u>Lemma 4.2</u>
*For the following special case of the discrete Gronwall inequality*

$$y_i \le a + b\sum_{j=0}^{i-1} y_j,$$

*with the initial condition* $y_0 = 0$*, the following inequality holds for* $y_i$*:*

$$(\forall i): y_i \le a(1+b)^{i-1}.$$

This lemma is proved by induction, and we omit the details. □

Continuing the proof of Theorem 4.1, we observe that, by Lemma 4.2, we have for the solution of (4.38):

$$(\forall i): |\varepsilon_i| \le h\left(1 + D^{(N)}h\right)^{\frac{t_i}{h}-1} \le F^{(N)}h \exp(TD^{(N)}) = O(h), \tag{4.40}$$

which completes the proof of Theorem 4.1. □

## 5. FINITE DIFFERENCE SCHEME FOR THE INFINITE-DIMENSIONAL CASE

We take a partition of the interval $[0,T]$, as we described in section 3. We approximate the values $\{x(t_i)\}_{i=0}^{M}$ of the solution of equation (1.6) on the partition $\{t_i\}_{i=0}^{M}$ with the solution $\{x_i\}_{i=0}^{M}$ of the following finite-difference equation:

$$x_i = x_0(t_i) + \sum_{n=1}^{\infty} \frac{h^n}{n!} \sum_{I_n(i)} f_n(t_i; t_{j_1}, \dots, t_{j_n}; x_{j_1}, \dots, x_{j_n}), \tag{5.1}$$

where $I_n(i)$ is defined in (4.2) .
We shall prove:

<u>Theorem 5.1</u>
*We assume that the properties* (i), (ii), *and* (iii), *stated in the condition of the Theorem 4.1 hold for the kernels* $f_n$ *in the equation (1.6). In addition, we require that for all* $N$ :

(iv) $B^{(N)} \le \overline{B}$ ; (v) $C^{(N)} \le \overline{C}$ ; (vi) $D^{(N)} \le \overline{D}$ ; (vii) $E^{(N)} \le \overline{E}$ ;

*where the constants* $\overline{B}$ , $\overline{C}$ , $\overline{D}$ , *and* $\overline{E}$ *are all finite.*
*Then for all* $i$ *the absolute value of the error* $|\varepsilon_i|$ *of the numerical approximation of the solution of equation (1.6) is of the order of* $h$ .

First, we prove two Lemmas:

<u>Lemma 5.1</u>
*Let* $x_i^{(N)}$ *be the approximate solution of the N-dimensional Volterra equation (1.5) by the approximation scheme of section 4, and* $x_i$ *be the approximate solution of the infinite-dimensional Volterra equation (1.6) by the equation (5.1).*
*Then for all* $i$ *we have* $|x_i^{(N)} - x_i| \to 0$, *as* $N \to \infty$.
<u>Proof</u>
Using definitions (4.1) and (5.1), we write

$$|x_i - x_i^{(N)}| \le \sum_{n=1}^{N} \frac{h^n}{n!} \sum_{I_n(i)} |f_n(t_i; t_{j_1}, \dots, t_{j_n}; x_{j_1}, \dots, x_{j_n}) - f_n(t_i; t_{j_1}, \dots, t_{j_n}; x_{j_1}^{(N)}, \dots, x_{j_n}^{(N)})| +$$
$$+ \sum_{n=N+1}^{\infty} \frac{h^n}{n!} \sum_{I_n(i)} |f_n(t_i; t_{j_1}, \dots, t_{j_n}; x_{j_1}, \dots, x_{j_n})| . \tag{5.2}$$

For the second term on the right hand side of (5.2) we obtain

$$c_N := \sum_{n=N+1}^{\infty} \frac{h^n}{n!} \sum_{I_n(i)} |f_n(t_i; t_{j_1}, \dots, t_{j_n}; x_{j_1}, \dots, x_{j_n})| \le \sum_{n=N+1}^{\infty} \frac{C_n T^n}{n!} \to 0, \tag{5.3}$$

as $N \to \infty$, by condition (v).
Next, for the first term in the right hand side of (5.2) we have, by applying condition (ii):

$$\sum_{n=1}^{N} \frac{h^n}{n!} \sum_{I_n(i)} | f_n(t_i; t_{j_1},...,t_{j_n}; x_{j_1},...,x_{j_n}) - f_n(t_i; t_{j_1},...,t_{j_n}; x_{j_1}^{(N)},...,x_{j_n}^{(N)}) | \leq$$

$$\leq \sum_{n=1}^{N} \frac{h^n}{n!} \sum_{I_n(i)} L_n \sum_{k=1}^{n} | x_{j_k} - x_{j_k}^{(N)} | \leq h \left[ \sum_{0 \leq j \leq i-1} | x_j - x_j^{(N)} | \right] D^{(N)} \tag{5.4}$$

Using (5.3) and (5.4) in (5.2), we obtain

$$| x_i - x_i^{(N)} | \leq c_N + hD^{(N)} \sum_{j=1}^{i-1} | x_j - x_j^{(N)} | . \tag{5.5}$$

By applying the discrete Gronwall inequality of Lemma 4.2 to (5.5), we have for each $N$ and for each $i$:

$$| x_i - x_i^{(N)} | < c_N e^{D^{(N)}T} . \tag{5.6}$$

Then for each $i$ we obtain from (5.6), using (5.3) and condition (vi):

$$\lim_{N\to\infty} | x_i - x_i^{(N)} | \leq \lim_{N\to\infty} (c_N e^{D^{(N)}T}) \leq e^{\overline{D}T} \lim_{N\to\infty} c_N = 0 , \tag{5.7}$$

as required. □

Lemma 5.2
*Let* $x(t)$ *be the exact solution of the infinite-dimensional Volterra equation (1.6), and* $x^{(N)}(t)$ *be the exact solution of the N-dimensional Volterra equation (1.5). Then under conditions* (i) *through* (vii) *of Theorem 5.1 we have for all* $t \in [0,T]$ $| x^{(N)}(t) - x(t) | \to 0$ *as* $N \to \infty$*, and convergence is uniform on the interval* $[0,T]$.

Proof
Using equations (1.5) and (1.6), we write

$$| x^{(N)}(t) - x(t) | \leq \sum_{n=1}^{N} \frac{1}{n!} \underbrace{\int_0^t ... \int_0^t}_{n \ times} | f_n(t, s_1,...,s_n, x^{(N)}(s_1),...,x^{(N)}(s_n)) -$$

$$- f_n(t, s_1,...,s_n, x(s_1),...,x(s_n)) | \, ds_n ... ds_1 +$$

$$+\sum_{n=N+1}^{\infty}\frac{1}{n!}\underbrace{\int_0^t\dots\int_0^t}_{n\ times}|f_n(t,s_1,\dots,s_n,x(s_1),\dots,x(s_n))|\,ds_n\cdots ds_1. \tag{5.8}$$

For the last term in the right hand side of inequality (5.8), we have, using condition (iii):

$$c^{(N)}:=\sum_{n=N+1}^{\infty}\frac{1}{n!}\underbrace{\int_0^t\dots\int_0^t}_{n\ times}|f_n(t,s_1,\dots,s_n,x(s_1),\dots,x(s_n))|\,ds_n\cdots ds_1\le\sum_{n=N+1}^{\infty}\frac{C_nT^n}{n!}\to 0 \tag{5.9}$$

as $N\to\infty$, by condition (v).
For the first term in the right hand side of inequality (5.8), we have, using condition (ii):

$$\sum_{n=1}^{N}\frac{1}{n!}\underbrace{\int_0^t\dots\int_0^t}_{n\ times}|f_n(t,s_1,\dots,s_n,x^{(N)}(s_1),\dots,x^{(N)}(s_n))-$$

$$-f_n(t,s_1,\dots,s_n,x(s_1),\dots,x(s_n))|\,ds_n\dots ds_1\le$$

$$\le\sum_{n=1}^{N}\frac{1}{n!}\underbrace{\int_0^t\dots\int_0^t}_{n\ times}L_n\sum_{k=1}^{n}\left|x^{(N)}(s_k)-x(s_k)\right|ds_n\dots ds_1\le D^{(N)}\int_0^t\left|x^{(N)}(s)-x(s)\right|ds \tag{5.10}$$

Next, using (5.9) and (5.10), we obtain for (5.8), by virtue of the continuous-time Gronwall inequality:

$$|x^{(N)}(t)-x(t)|\le c^{(N)}+D^{(N)}\int_0^t\left|x^{(N)}(s)-x(s)\right|ds. \tag{5.11}$$

We invoke the following result from [C1]:

*Suppose we have an inequality* $y(t)\le a+b\int_0^t y(s)ds$, *where* $y(s)\ge 0$ *for each* $s:0\le s\le t\le T$, *and* $a,b>0$. *Then the following inequality holds for each* $t\in[0,T]$*:*
$y(t)\le ae^{bt}$.

By applying the above result and condition (vi) to inequality (4.11), we have

$$0\le\lim_{N\to\infty}|x^{(N)}(t)-x(t)|\le e^{\bar{D}t}\lim_{N\to\infty}c^{(N)}. \tag{5.12}$$

Thus for all $t\in[0,T]$ we have from (5.12), using (5.9):

$$\lim_{N\to\infty} | x^{(N)}(t) - x(t) | = 0 , \qquad (5.13)$$

so that for each $t \in [0,T]$ $x^{(N)}(t) \to x(t)$ as $N \to \infty$. We observe that in fact this convergence is uniform on $[0,T]$, because we have for all $t \in [0,T]$:

$$0 \le | x^{(N)}(t) - x(t) | \le c^{(N)} e^{D^{(N)}T} \le c^{(N)} e^{\overline{D}T} . \qquad (5.14)$$

This completes proof of Lemma 5.2. □

Proof of Theorem 5.1

We have for each $N$:

$$| \varepsilon_i | = | x(t_i) - x_i | \le | x(t_i) - x^{(N)}(t_i) | + | x^{(N)}(t_i) - x_i^{(N)} | + | x_i^{(N)} - x_i | . \qquad (5.15)$$

Hence

$$| \varepsilon_i | \le \lim_{N\to\infty} | x(t_i) - x^{(N)}(t_i) | + \lim_{N\to\infty} | x^{(N)}(t_i) - x_i^{(N)} | + \lim_{N\to\infty} | x_i^{(N)} - x_i | . \qquad (5.16)$$

Since $t_i \in [0,T]$, we have, by Lemma 5.2:

$$\lim_{N\to\infty} | x(t_i) - x^{(N)}(t_i) | = 0 . \qquad (5.17)$$

Also, we have by Lemma 5.1:

$$\lim_{N\to\infty} | x_i^{(N)} - x_i | = 0 . \qquad (5.18)$$

Using (5.17) and (5.18), we obtain from (5.16):

$$| \varepsilon_i | \le \lim_{N\to\infty} | x^{(N)}(t_i) - x_i^{(N)} | = \lim_{N\to\infty} | \varepsilon_i^{(N)} | . \qquad (5.19)$$

We have from (4.40) for each $N$:

$$| \varepsilon_i^{(N)} | \le hF^{(N)} \exp(TD^{(N)}) \le h\left( \frac{\overline{D}T(M_0 + \overline{B} + \overline{C})}{2} + \overline{E} \right) \exp(T\overline{D}) , \qquad (5.20)$$

where we used conditions (iv)-(vii) of the Theorem. Using (5.19), we obtain

$$|\varepsilon_i| \le h\left(\frac{T\overline{D}}{2}(\overline{B}+\overline{C}+M_0)+\overline{E}\right)\exp(T\overline{D})\,, \tag{5.21}$$

and consequently $|\varepsilon_i| = O(h)$. □